\documentclass{elsart}
\usepackage{amsmath}
\usepackage[usenames]{color}
\usepackage{graphicx}

\usepackage{amssymb}

\begin{document}

\begin{frontmatter}



\title
{A number theoretical observation of  a 
 resonant interaction of Rossby waves
}


\author{Nobu Kishimoto}
\address{Research Institute for Mathematical Sciences, Kyoto University, Kyoto, 606-8502, Japan
}
\ead{nobu@kurims.kyoto-u.ac.jp}
\author{Tsuyoshi Yoneda
}
\address{Department of Mathematics, Tokyo Institute of Technology, Meguro-ku, Tokyo 152-8551, Japan 
}
\ead{yoneda@math.titech.ac.jp
}
\begin{abstract}
Rossby waves are generally expected to dominate the $\beta$ plane dynamics in geophysics, 
and here in this paper we give a number theoretical observation of  
the resonant interaction with a Diophantine equation.
The set of resonant frequencies does not have any frequency on the horizontal axis.
We also give several clusters of resonant frequencies.  
\end{abstract}


\end{frontmatter}


\section{Introduction}
\noindent
We consider three-wave interactions of the Rossby waves in a number theoretical approach.
Such waves are  observed in an incompressible two-dimensional flow on a $\beta$ plane (in geophysics). 
The $\beta$-plane approximation was first introduced by meteorologists 
(see \cite{C1,C2}) as a tangent plane of a sphere to approximately 
describe fluid motion on 
a rotating sphere, assuming that the Colioris parameter is a linear function 
of the latitude. 
A formal
derivation of the $\beta$-plane approximation is given in \cite{Pe}. 
It has been known that in
the incompressible two-dimensional flow on a $\beta$ plane, as time goes on, a zonal pattern emerges, consisting of alternating 
eastward and westward zonal flows, similar to the zonal band structure observed 
on Jupiter. 
From a physical point of view, one of the most important 
properties of the flow on a $\beta$ plane
linear waves 
called ``Rossby waves". The Rossby waves originate from 
the following dispersion relation (see \cite{YY} for example),
\begin{equation}
\omega=-\frac{\beta k_1}{k_1^2+k_2^2},
\end{equation}
where $\omega$ and $(k_1,k_2)$ are the angular frequency and the wavenumber vector. 
The Rossby waves
have been considered 
to play important roles in the dynamics of geophysical fluids (see \cite{Smith} for example).  
 In \cite{YY}, 
they proved a mathematical rigorous  theorem which supports the importance of the resonant
pairs of Rossby waves.  However, none of studies tried to consider such resonant waves in number theoretical approach,
and in this paper we attempt to consider it in an elementary number theory.
Let us be more precise. 
We define the wavenumber set consisting of wavenumbers in non-trivial
resonance as follows: 
\begin{defn}\label{wavenumber set} (Wavenumber set of non-trivial resonance.)\  
Let $\Lambda$ be a wavenumber set such that 
\begin{multline*}
\Lambda:=\bigg\{n\in \mathbb{Z}^2\quad\text{with}\quad n_1\not=0:\\
\frac{n_1}{n_1^2+n_2^2}-\frac{x}{x^2+y^2}-\frac{n_1-x}{(n_1-x)^2+(n_2-y)^2}=0,\\
 \text{for some}\ (x,y)\in\mathbb{Z}^2\ \text{with}\ x\not=0\ \text{and}\ n_1-x\not=0\bigg\}.
\end{multline*}
\end{defn}
The role of the above non-trivial resonance can be found in \cite{YY} in PDE sense.
Thus we omit to explain how it works to the two-dimensional flow on a  $\beta$ plane (in PDE  sense).
We would like to figure out the exact elements of $\Lambda$ without any numerical computation.
The following remark ensures that $\Lambda$ has at least infinite elements.
\begin{rem} (Infinite elements.)\ 
At least, 
$n=(n_1,n_2)=(m^4,m\ell^3)$ ($m,\ell\in\mathbb{N}$, $m\not=\ell$) is in $\Lambda$. In this case, we just take $(x,y)=(\ell^4,-m^3\ell)$.
Thus $\Lambda$ has at least infinite elements.
\end{rem}
\noindent
 $\Lambda$ itself is not only mathematically but also physically interesting.
In a turbulent flow, every wavenumber
component should have  nonzero energy. Suppose that the initial energy
distribution in a wavenumber space is isotropic. Two-dimensional
turbulence is known to transfer the energy from small to largescale
motions (energy inverse cascade). If there is no effect of rotation (no Coriolis effect), then the energy therefore becomes concentrated
isotropically around the origin in wavenumber space. However, if the rotation effect (Coriolis effect) is dominant, the energy transfer becomes governed by the resonant
interaction of Rossby waves $\Lambda$, and the number of resonant
triads gives a rough estimate of the strength of the nonlinear energy
transfer. Therefore, roughly speaking (in a physical point of view), the wavenumbers not in $\Lambda$ 
 are then expected to gain less energy compared
with wavenumbers in $\Lambda$.
In a numerical computation (see \cite{YY}), we can expect that $\Lambda$ has anisotropic distribution. 
Thus our aim  is  to know $\Lambda$ rigorously, and prove (in a number theoretical approach) that  
its distribution is anisotropic 
(however, it seems so difficult that we need to progress little by little).
For the first step,  in this paper, we give 
nonexistence of three wave interaction on $n_1$-axis by using a Diophantine equation. 
In Appendix, we give  suitable definitions to describe resonant points $\Lambda$, and give several specific resonant points. Up to now, the points were found one by one (not theoretically).
The main theorem is as follows:
\begin{thm}\label{nonexistence}
(Nonexistence of the three wave interaction on $n_1$-axis.)\ 
If $n_1,x,y\in\mathbb{Z}$ and
\begin{equation}\label{threewave}
\frac{n_1}{n_1^2}=\frac{x}{x^2+y^2}+\frac{n_1-x}{(n_1-x)^2+y^2},
\end{equation}
then $n_1x(n_1-x)=0$.
\end{thm}
\begin{rem}
In order to consider more general setting, namely, to figure out whether $(n_1,n_2)$ ($n_1,n_2\in\mathbb{Z}$, $n_1\not=0$) belongs to $\Lambda$ or not, we need to consider the following equality (just derived from Definition \ref{wavenumber set}):
\begin{multline*}
y^4-2n_2y^3-2x(n_1-x)y^2\\
+2n_2x\left(n_1-x+\frac{n_2^2}{n_1}\right)y-x(n_1-x)(x^2-n_1x+n_1^2+2n_2^2)-\frac{n_2^4x}{n_1}=0
\end{multline*}
for $x,y\in\mathbb{Z}$ with $x\not=0$.

\noindent
This equality might be related to ``elliptic curve"  more or less.
In this point of view, the ideas of  Mordell's theorem and  ``infinite descent" might be useful.
\end{rem}

\section{Proof of the Theorem \ref{nonexistence}.}
\noindent
Assume there is $n_1$ and $x$ such that $n_1x(n_1-x)\not=0$. Since $\Lambda$ is symmetric,  we can assume $n_1>x>0$.
From \eqref{threewave}, we see 
\begin{eqnarray*}
& & (x^2+y^2)\{(n_1-x)^2+y^2\}=n_1x\{(n_1-x)^2+y^2\}+n_1(n_1-x)(x^2+y^2)\\
&\Leftrightarrow&
y^4+\{x^2+(n_1-x)^2-n_1x-n_1(n_1-x)\}y^2\\
& &\ \ \ \ \ +\{x^2(n_1-x)^2-n_1x(n_1-x)^2-n_1x^2(n_1-x)\}=0\\
&\Leftrightarrow&
y^4-2x(n_1-x)y^2+x^2(n_1-x)^2-n_1^2x(n_1-x)=0\\
&\Leftrightarrow&
y^2=x(n_1-x)\pm\sqrt{x(n_1-x)}.
\end{eqnarray*}
Clearly, we do not treat complex numbers in this consideration, thus $n_1-x>0$. 
By $0<x(n_1-x)<n_1^2$ (we have already assumed that $n_1>x>0$, thus $n_1-x\leq n_1$),
we have $0<x(n_1-x)<n_1\sqrt{x(n_1-x)}$. Thus
\begin{equation*}
y^2=x(n_1-x)+n\sqrt{x(n_1-x)}.
\end{equation*}
Otherwise, $y$ becomes a complex number.
In particular, $x(n_1-x)=:p^2$ $(p\in\mathbb{N})$ and $p^2+n_1p$ are square numbers (if $x(n_1-x)$ is not square number,
then $y^2$ is not in $\mathbb{Z}$ and it is in contradiction).
Here, we can assume $x$ and $n_1$ are relatively prime. In fact, if the greatest common divisor is $d>1$,
we set $x'=x/d\in\mathbb{N}$ and $n_1'=n_1/d\in\mathbb{N}$ and then
\begin{equation*}
(y/d)^2=x'(n_1'-x')+n'\sqrt{x'(n_1'-x')}.
\end{equation*}
Since the left hand side of the  above equality is a rational number,  then $x'(n_1'-x')$ is a square number, namely, the right hand side 
is a natural number: $y':=y/d\in\mathbb{N}$.
Therefore we can regard  $n_1'$, $x'$ and $y'$ as $n_1$, $x$ and $y$. 
Since $x$ and $x(n_1-x)$ are relatively prime and $x(n_1-x)$ is a square number, $x$ and $n_1-x$ are also square numbers.
In fact, if either $x$ or $n_1-x$ is not square number, then (at least) two $p_j$ in the following expression 
\begin{equation*}
x(n_1-x)=p^2=p_1^2p_2^2\cdots p_N^2
\end{equation*}
($p_1,\cdots,p_N$ are prime numbers, and some $p_i$ and $p_j$ ($i\not=j$) may be the same)
must belong to both  $x$ and $(n_1-x)$.
In this case, $x$ and $n_1-x$ are not relatively prime.
Therefore
\begin{equation}\label{relatively prime}
x=q^2,\quad n_1=q^2+r^2,\quad p=qr,\quad q,r\in\mathbb{N}\ \text{are relatively prime}.
\end{equation}
We see that $q$, $r$ and $q^2+qr+r^2$ are all relatively prime.
For example, if $q$ and $q^2+qr+r^2$ are not relatively prime, there is a prime number $p$ such that 
$q=s_1p$ and  $q^2+qr+r^2=s_2p$ ($s_1,s_2\in\mathbb{N}$). 
Since $q^2+qr$ is multiple of $q$ (namely, multiple of $p$) then $r^2$ is also multiple of $p$.
However, if $r^2$ is multiple of $p$, then $r$ itself must be multiple of $p$.
This means that $q$ and $r$ are not relatively prime. It is in contradiction  to \eqref{relatively prime}.
Recall that $p^2+n_1p=q^2r^2+(q^2+r^2)qr=pr(q^2+qr+r^2)$ is a square number. Since  $q$, $r$ and $q^2+qr+r^2$ are all relatively prime,
we can rewrite
\begin{equation*}
q=s^2,\quad r=t^2,\quad s^4+s^2t^2+t^4=u^2,\quad s,t,u\in \mathbb{N}.
\end{equation*}
However it is in contradiction to the following lemma.
\begin{lem} (\cite{D})\ 
The following Diophantine equation
\begin{equation*}
X^4+X^2Y^2+Y^4=Z^2,\quad X,Y,Z\in\mathbb{Z} 
\end{equation*}
only have a trivial integer solution: $X=0$ or $Y=0$.
\end{lem}

\section{Appendix}
\noindent
In this section  we give several specific  resonant points. The points were found one by one (not theoretically).
In order to state the resonant points, ``cluster" concept is very useful.
Note that on a sphere case,  Kartashova and L'vov \cite{KL}  have already tried to classify several Rossby waves into clusters.
Let $\{\Omega^{finite}_j\}_j\subset\Lambda$ be a family of clusters composed by finite elements and 
$\{\Omega^{infty}_j\}_j\subset\Lambda$ be a family of clusters composed by infinite elements defined as follows:  
\begin{defn} (Clusters with finite elements)\  
Let $\{\Omega^{finite}_j\}_j\subset\Lambda$ be a family of wavenumber clusters  satisfying the following properties:
 
\begin{itemize}

\item

For any $n\in\Omega^{finite}_j$, there is  $(x,y)\in\Omega^{finite}_j$ 
such that $n$ and $(x,y)$ satisfy the definition of $\Lambda$.

\item

For any $n\not\in\Omega^{finite}_j$, there is no  $(x,y)\in\Omega^{finite}_j$ 
such that $n$ and $(x,y)$ satisfy the definition of $\Lambda$.

\item 

For each $j$, number of elements in $\Omega^{finite}_j$ is always finite.



\item

We set $\lambda_{1,j}:=\inf\{ |n|^2:n\in \Omega^{finite}_j\}$. 
Then $\lambda_{1,j}\leq \lambda_{1,k}$ ($j< k$).
If there is $j$ and $k$ $(j<k)$ such that $\lambda_{1,j}=\lambda_{1,k}$, again, we set $\lambda_{2,j}:=\inf\{|n|^2:n\in\Omega_j^{finite},\ |n|^2\not=\lambda_{1,j}\}$ and then $\lambda_{2,j}\leq \lambda_{2,k}$ ($j< k$). If  there is $j$ and $k$ ($j<k$) such that $\lambda_{2,j}=\lambda_{2,k}$, again,  we proceed the same manner.

\end{itemize}
\end{defn}
\begin{defn} (Clusters with infinite elements)\  
Let $\{\Omega^{infty}_j\}_j\subset\Lambda$ be a family  of wavenumber clusters satisfying the following properties:
 
\begin{itemize}

\item

For any $n\in\Omega^{infty}_j$, there is  $(x,y)\in\Omega^{infty}_j$ such that  $n$ and $(x,y)$ 
satisfy the definition of $\Lambda$.

\item

For any $n\not\in\Omega^{infty}_j$, there is no  $(x,y)\in\Omega^{infty}_j$ 
such that $n$ and $(x,y)$ satisfy the definition of $\Lambda$.

\item 

For each $j$, number of elements in $\Omega^{infty}_j$ is always infinity.

\item

We set $\lambda_{1,j}:=\inf\{ |n|^2:n\in \Omega^{infty}_j\}$. 
Then $\lambda_{1,j}\leq \lambda_{1,k}$ ($j< k$).
If there is $j$ and $k$ $(j<k)$ such that $\lambda_{1,j}=\lambda_{1,k}$, again, we set $\lambda_{2,j}:=\inf\{|n|^2:n\in\Omega_j^{infty},\ |n|^2\not=\lambda_{1,j}\}$ and then $\lambda_{2,j}\leq \lambda_{2,k}$ ($j< k$). If  there is $j$ and $k$ ($j<k$) such that $\lambda_{2,j}=\lambda_{2,k}$, again, we proceed the same manner.


\end{itemize}
\end{defn}
\begin{rem}\ 
\begin{itemize}

\item

We see $\Omega^{finite}_j\cap\Omega^{finite}_k=\emptyset\quad (j\not=k)$,
$\Omega^{infty}_j\cap\Omega^{infty}_k=\emptyset\quad (j\not=k)$ and $\Omega^{finite}_j\cap\Omega^{infty}_k=\emptyset$
\item

Theoretically,  each $\Omega^{finite}_j$ and $\Omega^{infty}_j$ ($j=1,2,\cdots$) are uniquely determined.

\item

We see
 
\begin{equation*}
\Lambda=\left(\cup_j\Omega_j^{finite}\right)\bigcup\left(\cup_j\Omega_j^{infty}\right).
\end{equation*} 
\end{itemize}
\end{rem}
In order to find  specific clusters, the following observation is useful.
For fixed $n$, 
we only need to see finite combinations of $(x,y)$ satisfying  the following inequality:
\begin{multline*}
\frac{|n_1|}{|n|^2} \le 
\frac{1}{\sqrt{x^2+y^2}} + \frac{1}{\sqrt{(n_1-x)^2+(n_2-y)^2}} \leq \\
\frac{2}{\min (\sqrt{x^2+y^2}, \sqrt{(n_1-x)^2+(n_2-y)^2})}.
\end{multline*}
Note that too large $x$ or $y$ will break  the above inequality.
Now we give several clusters: 
\begin{eqnarray*}
\Omega_1^{finite}&=&\{(1,11) ,(8,-34), (-9,23)\}\\
\Omega_2^{finite}&=&\{ (3,19), (32,-44), (-35,25) ,(8,26), (27,-51)\}\\
\end{eqnarray*}
and 
\begin{eqnarray*}
\Omega_1^{infty}&=&\{(1,-8) ,(3,-11), (5,25), (8,14), (13,13),\\
& &\ \ \ \ \ \ \  (15,10), (-16,-2), (27,-21), (-32,-4)\cdots\}.\\
\end{eqnarray*}
Since $(-32,-4)$ is two times $(-16,-2)$, we easily see that the above $\Omega_1^{infty}$ includes
\begin{eqnarray*}
& &\bigcup_{j=1}^\infty\{(j,-8j) ,(3j,-11j), (5j,25j), (8j,14j), (13j,13j),\\
& &\ \ \ \ \ \ \ \ \ \ \ 
 (15j,10j), (-16j,-2j), (27j,-21j), (-32j,-4j)\}.
\end{eqnarray*}
This means that $\Omega^{infty}_1$ has infinite elements. 
However, we could not figure out the exact elements in  $\Omega^{infty}_1$ so far.

{\it Acknowledgements.}\ 
The first author was partially supported by 
Grant-in-Aid for Young Scientists (B), No.~24740086,
Japan Society for the Promotion of Science.
The second author was partially supported by 
Grant-in-Aid for Young Scientists (B), No.~25870004,
Japan Society for the Promotion of Science.

\end{document}